\documentclass[11pt]{article}
\usepackage{amsmath,amssymb,amsthm,amsfonts,amscd,amsbsy}
\usepackage{diagrams}
\author{Peter J. Eccles and Mark Grant\footnote{Mark Grant is supported by a grant from
the United Kingdom Engineering and Physical Sciences Research Council.} \\
University of Manchester}

\date{}

\title{Bordism Groups of Immersions and Classes Represented by Self-Intersections}

\newcommand{\imm}{\looparrowright}
\newcommand{\bdm}{\begin{displaymath}}
\newcommand{\edm}{\end{displaymath}}
\newcommand{\R}{\mathbb{R}}
\newcommand{\Z}{\mathbb{Z}}
\newcommand{\I}{\mathcal{I}}
\newcommand{\co}{\colon\thinspace}
\newtheorem{thmH}{Theorem (Herbert)}[section]
\newtheorem{thm}[thmH]{Theorem}
\newtheorem{Cthm}[thmH]{Compression Theorem}
\newtheorem{prop}[thmH]{Proposition}
\newtheorem{propK}[thmH]{Proposition (Koschorke and Sanderson)}
\newtheorem{def.}[thmH]{Definition}
\newtheorem{cor.}[thmH]{Corollary}
\begin{document}
\maketitle
\begin{abstract}
We prove a geometrical version of Herbert's theorem \cite{He} by considering the self-intersection immersions of a self-transverse immersion up to bordism. This generalises Herbert's theorem to additional cohomology theories and gives a commutative diagram in the homotopy of Thom complexes. The proof is based on Herbert's and uses Koschorke and Sanderson's operations \cite{KS} and the fact that bordism of immersions gives a functor on the category of smooth manifolds and immersions.
\end{abstract}
\section{Introduction}

Given a self-transverse immersion $f\co M^{n-k}\imm N^{n}$ we define for each integer $r\geq 1$ the \emph{$r$-fold self intersection sets} of $f$, in $M$ and $N$ respectively, as
\begin{eqnarray*}
N_{r}&:=&\{ n\in N \mid |f^{-1}(n)|\geq r\}\subseteq N,\\
M_{r}&:=&f^{-1}(N_{r})\subseteq M.
\end{eqnarray*}
(Here and elsewhere the superscript denotes dimension, and all manifolds and maps are assumed smooth, so that for example $N^{n}$ is a smooth $n$-manifold while $N^{(n)}$ is the cartesian product of $n$ copies of $N$.) Using self-transversality of $f$ it can be shown (see \S 4.1) that each of these sets is the image of an immersed $(n-rk)$-manifold; that is, we have immersions
\bdm
\psi_{r}(f)\co\Delta_{r}(f)^{n-rk} \imm N^{n}, 
\edm
\bdm
\mu_{r}(f)\co\widetilde{\Delta}_{r}(f)^{n-rk} \imm M^{n-k},
\edm
whose images are $N_{r}$ and $M_{r}$ respectively. A general immersion $f\co M\imm N$ \emph{represents} a class in $H_{\ast}(N;\Z_{2})$ (and hence also in $H^{\ast}(N;\Z_{2})$ by Poincar\'{e} duality) by evaluating $f_{\ast}$ on the fundamental class of $M$. It is natural to ask how the homology classes represented by the above {\em self-intersection immersions} are related, for different values of $r$.

In his 1981 thesis \cite{He}, R.J. Herbert answered this question in the following satisfying way. Let $n_{r}\in H^{rk}(N;\Z_2)$ be the Poincar\'{e} dual of the class represented by $\psi_{r}(f)$, and let $m_{r}\in H^{(r-1)k}(M;\Z_2)$ be dual to the class represented by $\mu_{r}(f)$. Let $e\in H^{k}(M;\Z_2)$ be the Euler class of the normal bundle of $f$.
\begin{thmH}
\label{Herb} For $r\geq 1$,
\bdm
f^{\ast}n_{r} =  m_{r+1}+e\cup m_{r} \in H^{rk}(M;\Z_2).
\edm

\end{thmH}
Herbert's formula also holds in $H^{\ast}(M;\Z)$ when $M$ and $N$ are oriented and the codimension $k$ is even. Our aim here is to generalise this result by showing it to be true in any homology theory for which the constituent manifolds and normal bundles have an orientation, and to elucidate the simple geometric ideas behind Herbert's proof. To do this we study bordism groups of immersions as contravariant functors on the category of smooth manifolds and immersions. This turns out to be the ideal setting, as operations taking the bordism class of $f$ to that of $\psi_{r}(f)$ exist and are fairly well understood (see \cite{KS}, \cite{Vo}) and also the geometric meanings of constructions such as induced maps and products become apparent.

\section{Bordism of Immersions}

In this section we describe the functor given by bordism of immersions into a given manifold with a given structure on the normal bundle. The exposition follows closely that of Vogel \cite{Vo}. More details will appear in \cite{Gr}.
\subsection{Definitions}

Let $N^{n}$ be a connected manifold without boundary, and let $\zeta$ be vector bundle of dimension $k$ over a paracompact base space $X$. We consider immersions of closed manifolds into $N$ with a $\zeta$-structure on their normal bundles; that is, triples $(M^{n-k},f,v)$ where $M$ is a compact manifold without boundary, $f\co M\imm N$ is an immersion, and $v\co \nu(f)\to\zeta$ is a bundle map from the normal bundle of $f$. We may put a bordism relation on such triples in the following way. Two such triples $(M_0,f_0,v_0)$ and $(M_1,f_1,v_1)$ are \emph{bordant} if there exists a triple $(W^{n-k+1},F,V)$, where $W$ is a manifold with boundary, $F\co W\imm N\times I$ where $I$ is the unit interval $I=[0,1]$, and $V\co \nu(F)\to\zeta$ is a bundle map, and the following conditions are satisfied:

\textbf{1)} $\partial W=M_0\sqcup M_1$ ;

\textbf{2)} $F$ is transversal to $\partial (N\times I)=N\times\partial I$ ;

\textbf{3)} $F(\partial W)\subseteq N\times\partial I$ and $F|_{M_i}=f_i\times i$, for $i=0,1$ ;

\textbf{4)} $\nu(F)|_{M_i}\cong\nu(f_i)$ and $V|_{\nu(f_i)}=v_i$ .\\

Here it is understood that `=' means `equals up to diffeomorphism'. The resulting set of bordism classes will be denoted by $\mathcal{I}(N;\zeta)$, and the class of a triple $(M,f,v)$ will usually be written $[f]$ for brevity. It is clear that the set $\I(N;\zeta)$ has a commutative monoid structure arising from the disjoint union of immersions, and that the empty immersion provides a zero element. The complete homotopy classification of these monoids may be described (after Wells \cite{We}) as follows.
\begin{thm}
\label{htpy}
When $k=dim(\zeta)\geq 1$,
\bdm
\I(N;\zeta)\cong [N_+,T\zeta]_{S}\cong[N_+,QT\zeta].
\edm
When $k=dim(\zeta)=0$,
\bdm
\I(N;\zeta)\cong[ N_+,\bigvee^{\infty}_{r=1}(BS_r)_+].
\edm
\end{thm}
Here $T\zeta$ is the Thom space of $\zeta$, $N_+$ is the one-point compactification of $N$, and $QX$ is the direct limit of the spaces $\Sigma^{l}\Omega^{l}X$ (as $l\to\infty$) for pointed $X$. The space $BS_r$ is the classifying space of the symmetric group $S_r$.

In the case of positive codimension, Wells proved the above for $N$ the Euclidean space $\R^n$ and $\zeta =\gamma_k$, the universal $O(k)$-bundle. The result can be deduced for general $N$ and $\zeta$ from Thom's theorem on the bordism of embeddings (\cite{TH}) and the following useful theorem of differential topology (\cite{CT1}, \cite{CT3}).
\begin{Cthm}
\label{Cthm}
Let $g\co M^{n-k}\hookrightarrow N^{n}\times\R^l$ be an embedding with $l$ linearly independent normal vector fields, so that $\nu(g)\cong\nu\oplus\varepsilon^l$ for some $k$-dimensional vector bundle $\nu$. If $k\geq 1$ then $g$ is isotopic to an embedding $g'$ such that the composition
\bdm
M\stackrel{g'}{\hookrightarrow}N\times\R^l\stackrel{pr}{\rightarrow}N
\edm
is an immersion $f\co M\imm N$. Moreover, $\nu(f)\cong\nu$.
\end{Cthm}
 When the codimension is zero the Compression Theorem fails. However, the result follows from the observation that a codimension zero immersion of a closed manifold is nothing but a finite covering space. A direct proof covering both cases and using configuration space models has been given by Koschorke and Sanderson \cite{KS}.

Since $\I(N;\zeta)$ is isomorphic to a stable homotopy group when $k\geq1$, it has the structure of an Abelian group. Commutativity is clear, but it seems we cannot give a canonical representative for the inverse of the class of a triple $(M,f,v)$.

There is also an \emph{external product} pairing
\bdm
\I(N;\zeta)\times\I(N';\zeta')\stackrel{\times}{\to}\I(N\times N';\zeta\times\zeta'),
\edm
given fairly obviously by the cartesian product of representatives. This product corresponds to the smash product of stable homotopy classes under the isomorphism of Theorem \ref{htpy}.
\subsection{Functoriality}
Theorem \ref{htpy} shows that $\I(-;-)$ is a homotopy bifunctor, covariant in $\zeta$ and contravariant in $N$. This may be interpreted in terms of the differential topology, as we shall now explain. 

Let $\mathcal{C}_k$ be the category with objects the $k$-dimensional vector bundles and with morphisms the bundle maps which map fibres isomorphically. It is clear that a bundle map $\eta\co \zeta\to\xi$ induces a map $\eta_{\ast}\co \I(N;\zeta)\to\I(N;\xi)$ so that $\I(N;-)$ gives a covariant functor from $\mathcal{C}_k$ to the category of commutative monoids. Also note that bundle homotopic bundle maps induce the same map of monoids, so that we do indeed have a homotopy functor.

More subtle is the contravariance in $N$. Let $\mathcal{D}_0$ denote the category with objects the smooth manifolds without boundary and with morphisms the immersions. A homotopy in this category is a regular homotopy (a smooth map $F\co M\times I\to N$ such that at each stage $t\in I$ the map $F(-,t)\co M\to N$ is an immersion). 
\begin{prop}
\label{func}
$\I(-;\zeta)$ is a contravariant homotopy functor from $\mathcal{D}_0$ to the category of commutative monoids.
\end{prop}
\begin{proof}
Given a class $[f]$ in $\I(N;\zeta)$ and an immersion $g\co Q^{n-l}\imm N^{n}$, we may choose a representative $f'$ of $[f]$ which is transverse to $g$ and regularly homotopic to $f$; then the map $\delta$ in the following pullback diagram is an immersion.
\begin{diagram}[height=2.5em]
\Theta^{n-k-l} &\rTo^{\rho} & M^{n-k} \\
\dTo<{\delta} & &\dTo>{f'} \\
Q^{n-l} &\rTo^{g} & N^{n}
\end{diagram}
It can be shown (using Theorem \ref{Cthm}) that $\nu(\delta)\cong \rho^{\ast}\nu(f')\cong\rho^{\ast}\nu(f)$, and so we have a bundle map $\overline{\rho}\co \nu(\delta)\to\nu(f)$. This allows us to define the induced map 
\bdm
g^{\ast}\co \I(N;\zeta)\to\I(Q;\zeta)
\edm
by sending the class of $(M,f,v)$ to the class of $(\Theta,\delta,v\circ\overline{\rho})$. It is easily checked that this is a well-defined map of monoids, and functoriality is evident. One may also verify that pulling back by regularly homotopic immersions of $Q$ in $N$ results in bordant immersions in $Q$.
\end{proof}
\textbf{Remarks.}\\

 \textbf{i)} Let $g$ be an immersion and $\eta$ a bundle map. Then the induced maps $g^{\ast}$ and $\eta_{\ast}$ are seen to correspond under Theorem \ref{htpy} to composition with the stable homotopy classes of $g_+\co Q_{+}\to N_+$ and $T\eta\co T\zeta\to T\xi$.

\textbf{ii)} Functoriality and the diagonal embedding $\triangle\co N\to N\times N$ gives an \emph{internal product}
\bdm
\I(N;\zeta)\times\I(N;\zeta')\stackrel{\cup}{\to}\I(N;\zeta\times\zeta')
\edm
by setting $[f]\cup [f']=\triangle^{\ast}([f]\times [f'])$.
\subsection{Relative groups}
With a little care we may extend the definition of $\I(-;\zeta)$ to the category $\mathcal{D}$ of smooth manifolds with boundary and immersions.  \emph{The immersions in the domain category are not required to preserve boundaries}. Let $N$ have boundary $\partial N$. Now the data are triples $(M,f,v)$ where $M$ is a compact manifold with boundary $\partial M$, and $f\co M\imm N$ meets $\partial N$ transversally with $f^{-1}(\partial N)=\partial M$. The precise definition of a bordism between such triples involves the idea of a manifold with corners \cite{lau00}; the details will appear in \cite{Gr}. Call the resulting set of bordism classes $\I(N,\partial N;\zeta)$.

\section{Relations to Cohomology Theories}
In this section we shall outline the relation between the bifunctor $\I(-;-)$ and the various generalised (co)homology theories known as (co)bordism (further details may be found in the references). Such theories arise from the various \emph{Thom spectra} (see for example \cite{St}, \cite{Sw}) which are constructed as follows. 

Suppose we are given a family $\textbf{X} =\{X_{k},f_{k},g_{k}\}_{k\geq 0}$, where for each $k$ we have a space $X_{k}$, a fibration $f_{k}\co X_{k}\to BO(k)$ and a map $g_{k}\co X_{k}\to X_{k+1}$ making the following diagram commute.
\begin{diagram}[height=2.5em]
X_{k} & \rTo^{f_{k}} & BO(k) \\
\dTo<{g_{k}} & & \dTo>{Bi_{k}} \\
X_{k+1} &\rTo^{f_{k+1}} & BO(k+1)
\end{diagram}
Here $i_k\co O(k)\hookrightarrow O(k+1)$ is the inclusion homomorphism, and if $\gamma_k$ is the universal $O(k)$-bundle over $BO(k)$, we have $(Bi_{k})^{\ast}\gamma_{k+1}\cong\gamma_k\oplus\varepsilon^1$. From these data we get a sequence of bundles $\boldsymbol{\zeta} =\{\zeta_k :=f_{k}^{\ast}\gamma_k\}_{k\geq 0}$, with $dim(\zeta_k)=k$, and bundle maps $\overline{g}_k\co\zeta_k\oplus\varepsilon^1\to\zeta_{k+1}$. This gives rise to a \emph{Thom spectrum} $\textsf{M}\boldsymbol{\zeta}$, with $k$-th space $T\zeta_k$, and hence to generalised cohomology and homology theories $\textsf{M}\boldsymbol{\zeta}^{\ast}$ and $\textsf{M}\boldsymbol{\zeta}_{\ast}$. 

We suppose further that $\textsf{M}\boldsymbol{\zeta}$ is a ring spectrum, the multiplicative structure coming from compatible bundle maps $\zeta_k\times\zeta_l\to\zeta_{k+l}$ inducing maps $T\zeta_k\wedge T\zeta_l\to T\zeta_{k+l}$ for all $k,l\geq 0$. This furnishes our (co)bordism theory with cup and cap products. 

We shall need some notion of orientability of manifolds with respect to the spectrum $\textsf{M}\boldsymbol{\zeta}$.
\begin{def.}
A \emph{$\boldsymbol{\zeta}$-manifold} is a smooth manifold $N$ together with an equivalence class of $\boldsymbol{\zeta}$-structures on its stable normal bundle. Thus we have an embedding $g\co N^n \hookrightarrow\R^{n+l}$ and a homotopy class of bundle maps  $\nu(g)\to\zeta_l$, where $l$ is large.
\end{def.}

\subsection{Quillen's Geometric Cobordism and Duality}
Let $N^n$ be a closed manifold. It is well known that the bordism theory $\textsf{M}\boldsymbol{\zeta}_{\ast}$ has a geometric interpretation, given by Atiyah \cite{At}. Here $\textsf{M}\boldsymbol{\zeta}_{n-k}(N)$ is the set of bordism classes of singular maps $f\co M^{n-k}\to N^n$, where $M$ is a closed $\boldsymbol{\zeta}$-manifold. Less well known is Quillen's geometric interpretation of the cobordism group $\textsf{M}\boldsymbol{\zeta}^k (N)$, as equivalence classes of $\boldsymbol{\zeta}$-oriented codimension $k$ maps \cite{Qu}.
\begin{def.}
A map of closed manifolds $f\co M^{n-k}\to N^n$ is {\em $\boldsymbol{\zeta}$-orientable} if it has a factorisation
\bdm
M\stackrel{i}{\longrightarrow} E\stackrel{p}{\longrightarrow} N,
\edm
where $p\co E\to N$ is an $(r-k)$-dimensional smooth vector bundle over $N$ for some large $r$ (hence $E$ is a $(n-k+r)$-manifold), and $i$ is an embedding with a $\zeta_r$-structure on its normal bundle. A {\em $\boldsymbol{\zeta}$-orientation} of $f$ is a class of such factorisations with respect to a suitable equivalence relation.
\end{def.}

Quillen puts a cobordism relation on the class of $\boldsymbol{\zeta}$-oriented maps $f\co M\to N$ of codimension $k$. Then a generalisation of Thom's original proof for the coefficient groups shows that the resulting set of cobordism classes is isomorphic to 
\bdm
\textsf{M}\boldsymbol{\zeta}^{k}(N)\cong\lim_{r\to\infty}[ \Sigma^r N_+,T\zeta_{k+r}].
\edm
For the remainder of this section, we suppose $N$ is a closed $\boldsymbol{\zeta}$-manifold. This gives a Poincar\'{e}-Atiyah duality isomorphism $\textsf{M}\boldsymbol{\zeta}^{k}(N)\cong\textsf{M}\boldsymbol{\zeta}_{n-k}(N)$ as in \cite{At}, which is now described geometrically by the identity on representatives. A map $f\co M\to N$ between $\boldsymbol{\zeta}$-manifolds has a canonical $\boldsymbol{\zeta}$-orientation, and so can represent both a bordism and a cobordism class; these classes are dual to each other.
\subsection{$\I(N;\zeta_k)$ in this setting}
Now let an element $[f]\in\I(N;\zeta_k)$ be represented by the triple $(M^{n-k},f,v)$. Choose an embedding $f_1\co M\hookrightarrow\R^l$, for $l$ large. Let $i\co M\hookrightarrow N\times\R^l$ be the embedding given by $i(m)=(f(m),f_1(m))$ for $m\in M$. This gives a factorisation of $f$ as
\bdm
M\stackrel{i}{\longrightarrow} N\times\R^l \stackrel{pr}{\longrightarrow} N,
\edm
and the sequence of bundle maps
\bdm
\nu(i)\cong\nu(f)\oplus\varepsilon^l \stackrel{v\oplus 1}{\longrightarrow}\zeta_k\oplus\varepsilon^l \to\zeta_{k+l}
\edm
shows that we in fact have a $\boldsymbol{\zeta}$-orientation of $f$. After checking that a bordism between immersions with $\zeta_k$-structure gives a cobordism between $\boldsymbol{\zeta}$-oriented maps, we have a well defined map 
$\Theta_k\co \I(N;\zeta_k)\to\textsf{M}\boldsymbol{\zeta}^{k}(N)$, which is seen to coincide with the map
\bdm
[\Sigma^l N_+ ,\Sigma^l T\zeta_k]\to [\Sigma^l N_+ ,T\zeta_{k+l}]
\edm
given by the maps in the spectrum $\textsf{M}\boldsymbol{\zeta}$. 

In his paper \cite{Qu}, Quillen also gives the geometric constructions of induced maps, addition and products in the theory $\textsf{M}\boldsymbol{\zeta}^{\ast}$, and these are completely analogous to the constructions we gave in the case of immersions in the last section. Using this observation, or by examining the maps in homotopy theory, we see that each $\Theta_k$ is a homomorphism of monoids (groups for $k\geq 1$), and a natural transformation of functors on the category $\mathcal{D}$. The multiplicativity is expressed in the following diagram.
\begin{diagram}[height=2.5em]
\I(N;\zeta_k)\times\I(N;\zeta_l) &\rTo^{\cup} & \I(N;\zeta_{k+l}) \\
\dTo>{\Theta_k\times\Theta_l} & & \dTo>{\Theta_{k+l}} \\
\textsf{M}\boldsymbol{\zeta}^{k}(N)\times\textsf{M}\boldsymbol{\zeta}^{l}(N) &\rTo^{\cup} &\textsf{M}\boldsymbol{\zeta}^{k+l}(N)
\end{diagram}

There is also a homomorphism from $\I(N;\zeta_{k})$ to $\textsf{M}\boldsymbol{\zeta}_{n-k}(N)$. This is obtained by regarding an immersion $f\co M^{n-k}\imm N^n$ as a singular manifold in $N$. Note that $M$ is a $\boldsymbol{\zeta}$-manifold, since its stable normal bundle receives a $\boldsymbol{\zeta}$-structure from that of $\nu(f)$ and the stable normal bundle of $N$. The image of a class $[f]$ in $\textsf{M}\boldsymbol{\zeta}_{n-k}(N)$ is dual to its image in $\textsf{M}\boldsymbol{\zeta}^{k}(N)$, since both are represented by the map $f\co M\to N$. Finally, suppose we have a Thom class $t\co \textsf{M}\boldsymbol{\zeta}\to\textsf{A}$, that is, a ring map between ring spectra. Let $\textsf{A}^{\ast}$, $\textsf{A}_{\ast}$ denote the cohomology and homology theories arising from the spectrum $\textsf{A}$. Then $t$ gives a degree zero multiplicative natural transformation of theories which maps Thom class to Thom class. This also furnishes $N$ with an $\textsf{A}$-orientation, and hence we have Poincar\'{e} Duality with $\textsf{A}$ coefficients and the following commutative diagram.
\begin{diagram}[height=2.5em]
\mathcal{I}(N;\zeta_k) &\rTo^{\Theta_k} &\textsf{M}\boldsymbol{\zeta}^{k}(N)&\rTo^{\overline{t}} &\textsf{A}^{k}(N) \\
& \rdTo & \dTo>{\cong} & &\dTo>{\cong} \\
& &\textsf{M}\boldsymbol{\zeta}_{n-k}(N)&\rTo_{\underline{t}} &\textsf{A}_{n-k}(N)
\end{diagram}

{\bf Examples}

{\bf (i)} The $\textsf{MO}$ spectrum of unoriented bordism, and the so-called `universal' Thom class $t\co \textsf{MO}\to\textsf{H}(\Z_2)$.

{\bf (ii)} The $\textsf{MSO}$ spectrum of oriented bordism and the Thom class $t\co \textsf{MSO}\to\textsf{H}(\Z)$.

{\bf (iii)} The $\textsf{MU}$ spectrum of complex bordism and the Conner-Floyd map $t\co \textsf{MU}\to\textsf{K}$ to complex $K$-theory.

\section{Some Operations in the Bordism of Immersions}

Following Vogel \cite{Vo} and Koschorke and Sanderson \cite{KS}, we now introduce operations (natural transformations of set-valued co-functors)
\bdm
\psi^r\co \I(N;\zeta)\to\I(N;\mathcal{S}_r(\zeta )),
\edm
for each integer $r\geq 0$, where $\mathcal{S}_r(\zeta)$ is a bundle to be defined. These operations take the bordism class of an immersion $f\co M\imm N$ to the bordism class of its $r$-fold self-intersection immersion $\psi_r(f)\co \Delta_r(f)\imm N$.
\subsection{Constructions and Properties}
Let $f\co M^{n-k}\imm N^n$ be a self-transverse immersion with a $\zeta$-structure on its normal bundle. Write $f^{(r)}|$ for the restriction of $f^{(r)}$ to the open submanifold $F(M;r)\subset M^{(r)}$ of ordered $r$-tuples of {\em distinct} points in $M$. By the self-transversality of $f$, this immersion into $N^{(r)}$ is transverse to the diagonal embedding $\triangle\co N\hookrightarrow N^{(r)}$. We form the pullback square.
\begin{diagram}[height=2.5em]
\overline{\Delta}_r(f)^{n-rk} & \rTo & F(M;r) \\
\dTo & &\dTo>{f^{(r)}|} \\
N &\rTo^{\triangle} & N^{(r)}
\end{diagram}
The manifold
\bdm
\overline{\Delta}_r(f)^{n-rk}=\{(m_1,\ldots ,m_r)\in F(M;r)\mid f(m_1)=\ldots =f(m_r)\}
\edm
admits free actions of the symmetric groups $S_r$ and $S_{r-1}$, by permuting the last $r$ and $(r-1)$ co-ordinates respectively. Call the respective quotient manifolds $\Delta_{r}(f)^{n-rk}$ and $\widetilde{\Delta}_{r}(f)^{n-rk}$. We may now define the {\em $r$-fold self-intersection immersion of $f$ in $N$} as
\bdm
\psi_{r}(f)\co\Delta_{r}(f)^{n-rk} \imm N^{n} 
\edm
\bdm
[m_1,\ldots ,m_r]\mapsto f(m_1)=\ldots =f(m_r).
\edm
We define the {\em $r$-fold self-intersection immersion of $f$ in $M$} as
\bdm
\mu_{r}(f)\co\widetilde{\Delta}_{r}(f)^{n-rk} \imm M^{n-k}
\edm
\bdm
(m_1,[m_2,\ldots ,m_r])\mapsto m_1.
\edm
These are indeed immersions of compact manifolds (as shown in \cite{He}); we shall need to know the structure of their normal bundles.

Let $ES_r$ be a contractible space with a free action of $S_r$ on the right (for example, the total space of a universal principle $S_r$-bundle). For a vector bundle $p\co E(\zeta)\to B(\zeta)$, we define $\mathcal{S}_r(\zeta)$ to be the bundle
\bdm
1\times_{S_r}p^{(r)}\co ES_r\times_{S_r}E(\zeta)^{(r)}\longrightarrow ES_r\times_{S_r}B(\zeta)^{(r)}.
\edm
Here we quotient out by the diagonal action of $S_r$ on the product, where $S_r$ acts on $X^{(r)}$ on the left by permuting the factors. It is standard that there is a homeomorphism
\bdm
T\mathcal{S}_r(\zeta)\approx D_r T\zeta
\edm
where $D_r X$ is the $r$-adic construction (or equivariant half-smash product) on the pointed space $X$.

It may now be seen fairly easily that a $\zeta$-structure on the normal bundle of $f$ induces a $\mathcal{S}_r(\zeta)$ structure on the normal bundle of $\psi_r(f)$. The fibre of this normal bundle at a point $[m_1,\ldots ,m_r]\in\triangle_{r}(f)$ is the unordered direct sum of the fibres of the normal bundle of $f$ at the points $m_1,\ldots ,m_r$. (Similarly one sees that the normal bundle of $\mu_r(f)$ has a $\mathcal{S}_{r-1}(\zeta)$-structure). Hence we may define, for $r\geq 1$, an operation
\bdm
\psi_r\co \I(N;\zeta)\to\I(N;\mathcal{S}_r(\zeta ))
\edm
by setting $\psi_r [f]=[\psi_r(f')]$, where $f'$ is a self-transverse representative of $[f]$. This is well defined, since a bordism between self-transverse representatives may be chosen self-transverse, and this defines a bordism between the self-intersection immersions. 

Let $\star$ denote the trivial 0-dimensional bundle over a point, and define $\psi_0[f]$ to be $[\textbf{1}_N \co N\imm N]\in\I(N;\star)$. Note that this class acts as an identity for the $\cup$-product. The operations $\psi_r$ for $r\geq 0$ satisfy the following properties.
\begin{propK}
\label{properties}
Let $[f],[g]\in\I(N;\zeta)$.\\
{\bf (i)}(Naturality) If $h\co Q\imm N$, then
\bdm
h^{\ast}\psi_r[f]=\psi_r h^{\ast}[f]\in\I(Q;\mathcal{S}_r(\zeta )). 
\edm
{\bf (ii)} $\psi_1 [f]=[f]$, and $\psi_r [f]=0$ for $r>1$ if $[f]$ can be represented by an embedding.\\
{\bf (iii)} (Cartan formula)
\bdm
\psi_r([f]+[g])=\sum_{i=0}^{r} \psi_{r-i} [f]\cup\psi_i[g]\in \I(N;\mathcal{S}_r(\zeta )).
\edm
{\bf (iv)} If $f$ is self-transverse,
\bdm
\psi_r[\mu_2(f)]=[\mu_{r+1}(f)]\in\I(M;\mathcal{S}_r(\zeta )).
\edm
\end{propK}
\begin{proof}
Properties {\bf (i)},{\bf (ii)}, and {\bf (iv)} are direct consequences of the definitions. Property {\bf (iii)} is slightly less obvious, but the moral is that ``If $f\sqcup g$ is self-transverse then an $r$-fold self-intersection of $f\sqcup g$ is the intersection of an $(r-i)$-fold self-intersection of $f$ with an $i$-fold self-intersection of $g$". Notice we have made use of bundle maps $\mathcal{S}_{r-i}(\zeta)\times\mathcal{S}_i (\zeta)\to\mathcal{S}_{r}(\zeta)$, which exist thanks to the product $ES_{r-i}\times ES_i\to ES_{r}$ coming from concatenation of permutations, so that the formula ends up in $\I(N;\mathcal{S}_r(\zeta ))$.
\end{proof}
\subsection{The James-Hopf maps $h^r\co QX\to QD_r X$}

Recall that for a path-connected, pointed space $X$ there is a stable homotopy equivalence
\bdm
QX\simeq_{S}\bigvee_{r\geq 1}D_r X,
\edm
the so-called `Snaith splitting'.

This result rests on the existence of certain {\em James-Hopf maps}
\bdm
h^r\co QX\to QD_r X,\qquad r\geq 1,
\edm
which may be defined combinatorially using models for the space $QX$ (see for example \cite{bar74}). The following theorem, which has been proved by Vogel \cite{Vo} and by Sz\H{u}cs \cite{Sz1}, \cite{Sz2} (see also \cite{KS}), shows how these maps provide a homotopy interpretation of the operations $\psi_r$.
\begin{thm}
The following diagram commutes.
\begin{diagram}[height=2.5em]
\I(N;\zeta) &\rTo^{\cong} & [N_+ ,QT\zeta ] \\
\dTo<{\psi_r} & &\dTo>{(h^r)_*} \\
\I(N;\mathcal{S}_r(\zeta)) &\rTo^{\cong} & [N_+ ,QD_r T\zeta ]
\end{diagram}
\end{thm}
\section{Proving Herbert's Theorem}
We may now prove the analogue of Herbert's Theorem \ref{Herb} in $\I(M;\mathcal{S}_r(\zeta))$. Herbert's original proof in \cite{He} was encumbered by the fact that homology classes are represented by singular simplices rather than immersed manifolds. The relative simplicity of our proof stems from the fact that bordism of immersions is in some sense closer to the geometry. Our method is to first prove the double point case, then deduce the general case using the properties of $\psi_r$ listed in Proposition \ref{properties}.
\subsection{Proofs}
Let $f\co M\imm N$ be an immersion of closed manifolds with a $\zeta$-structure on its normal bundle. Therefore $f$ represents a class in $\I(N;\zeta)$ and $\mu_2(f)$ a class in $\I(M;\zeta)$. Let $D\nu(f)$ denote the closed disc bundle of the normal bundle of $f$, and $i\co M\hookrightarrow D\nu(f)$ its zero section.
\begin{prop}
\label{doub}
 \bdm
 f^{\ast}[f]=[\mu_2(f)]+e\in\I(M;\zeta),
 \edm
 where $e:=i^*[i]$ is the {\em Euler class} of $\nu(f)$.
\end{prop}
\begin{proof}
We may extend $f$ to an immersion $F\co D\nu(f)\imm N$ of the closed disc normal bundle of $f$ which is injective on the fibres. This immersion $F$ does not represent a bordism class, since it does not preserve boundaries. However, it is a morphism in the domain category $\mathcal{D}$ (see \S 2.3); and since $f=F\circ i$, we have
\bdm
f^{\ast}[f]=i^* F^*[f]
\edm
by functoriality. So we first give an explicit immersion which represents $F^*[f]\in\I(D\nu(f), \partial D\nu(f);\zeta )$.

Since $F$ has codimension zero, it is automatically transverse to $f$. The pull back construction of \ref{func} gives rise to a manifold
\bdm
\Theta =\{ (v,m)\in D\nu(f)\times M\mid F(v)=f(m)\},
\edm
and an immersion $\delta\co\Theta\imm D\nu(f)$ with $\delta (v,m)=v$. Having insisted that $F$ be injective on fibres, we can partition $\Theta$ as a disjoint union $\Theta=\Theta_0\sqcup\Theta_1$, where $\Theta_0 =\{(i(m) ,m)\}$, and 
\bdm
\Theta_1 =\{(v,m)\mid v\in\nu(f)_{m_1},\, m_1\neq m\}.
\edm
Clearly $\Theta_0 =M$ and the restriction of $\delta$ to $\Theta_0$ is identified as $i\co M\hookrightarrow D\nu(f)$. Putting $g:=\delta |_{\Theta_1}\co \Theta_1 \imm D\nu(f)$, we have
\bdm
F^*[f]=[i]+[g]\in\I(D\nu(f), \partial D\nu(f);\zeta ).
\edm

Applying $i^*$ now gives
\bdm
f^*[f]=e+i^*[g]
\edm
Now we may assume that $f$ is self-transverse. This implies that $i$ is transverse to $g$, making it easy to check that the pullback of $g$ by $i$ is indeed the immersion $\mu_2(f)\co \widetilde{\Delta}_{2}(f)\imm M$.
\end{proof}
\begin{thm}
\label{main}
For $r\geq 1$,
\bdm
f^* [\psi_r(f)]=[\mu_{r+1}(f)]+e\cup [\mu_r(f)]\in\I(M;\mathcal{S}_r(\zeta)).
\edm
\end{thm}
\begin{proof}
This follows from the double point case above using Proposition \ref{properties}. The argument is as follows:
\begin{eqnarray*}
f^*[\psi_r(f)] &=&\psi_r f^*[f] \quad\mbox{by {\bf (i)} of \ref{properties}}\\
& = &\psi_r ([\mu_2(f)]+e) \quad\mbox{by \ref{doub}}\\
& = &\sum_{i=0}^{r} \psi_{r-i} [\mu_2(f)]\cup\psi_i[e] \quad\mbox{by {\bf (iii)}}\\
& = & \psi_r[\mu_2(f)]+\psi_{r-1}[\mu_2(f)]\cup e \quad\mbox{by {\bf (ii)}}\\
& = & [\mu_r+1(f)]+e\cup [\mu_r(f)] \quad\mbox{by {\bf (iv)}}.
\end{eqnarray*}
\end{proof}
Using the diagram at the end of \S 3, immersions represent classes in certain cohomology theories according to their structure, and we may deduce Herbert's theorem in these theories. Let $\textsf{M}\boldsymbol{\zeta}$ be the Thom spectrum with $k$-th space $T\zeta_k$, and let $f\co M\imm N$ have a $\zeta_k$-structure on its normal bundle. We remark without proof that our Euler class $e:=i^*[i]$ maps to the cobordism Euler class of $\nu(f)$ under the homomorphism $\Theta_k\co\I(M;\zeta_k)\to\textsf{M}\boldsymbol{\zeta}^k (M)$ (compare the definition at the start of \S 2 of \cite{Qu}). Suppose there is a bundle map from $\mathcal{S}_r(\zeta_k)$ to $\zeta_{rk}$.
\begin{cor.}
With the above assumptions, Herbert's theorem holds in $\I(M;\zeta_{rk})$ and hence in $\textsf{M}\boldsymbol{\zeta}^{rk}$.
\end{cor.}
 {\bf Examples}
 
{\bf (i)} The $\textsf{MO}$ spectrum. Let $\gamma_k$ be the universal $O(k)$ bundle; then there is a classifying bundle map $\mathcal{S}_r(\gamma_k)\to\gamma_{rk}$ for every $r\geq 1$, and hence \ref{main} is true also in $\textsf{MO}^*$ and $\textsf{H}(\Z_2)^*$. This recovers Theorem \ref{Herb}.

{\bf (ii)} The $\textsf{MSO}$ spectrum. Let $\widetilde{\gamma_k}$ be the universal $SO(k)$ bundle. When $k$ is even, the bundle $\mathcal{S}_r(\widetilde{\gamma_k})$ is orientable, so admits a map to $\widetilde{\gamma_{rk}}$. Hence we deduce Theorem \ref{main} in the theories $\textsf{MSO}^*$ and $\textsf{H}(\Z )$ when the normal bundle of $f$ is orientable and of even dimension.

{\bf (iii)} The $\textsf{MU}$ spectrum. Let $\gamma^U_j$ be the universal $U(j)$ bundle (of real dimension $2j$). The bundle $\mathcal{S}_r(\gamma^U_j)$ has a complex structure constructed from that of $\gamma^U_j$. Hence if $f$ is an immersion of even codimension with a complex structure on its normal bundle, we obtain Theorem \ref{main} in the theories $\textsf{MU}^*$ and $\textsf{K}^*$.
\subsection{Homotopy of Thom Spaces}
By virtue of Theorem \ref{htpy}, bordism classes of immersions are represented by stable homotopy classes of maps into Thom spaces. We may therefore translate the geometric result Theorem \ref{main} into a homotopy commutative diagram of representing maps.
 Write $\underline{f}\co N_+\to QT\zeta$ for the stable adjoint of a map representing the class $[f]\in\I(N;\zeta )$. By \S 4.2 and the remarks following Proposition \ref{func}, the class $f^*\psi_r[f]$ is represented by the composition
\bdm
M_+\stackrel{f_+}{\longrightarrow}N_+\stackrel{\underline{f}}{\longrightarrow}QT\zeta\stackrel{h^r}{\longrightarrow}QD_r T\zeta.
\edm
Similarly, Proposition \ref{properties} {\bf (iv)} says that $\underline{\mu_{r+1}(f)}$ is the composition
\bdm
M_+\stackrel{\underline{\mu_2(f)}}{\longrightarrow}QT\zeta\stackrel{h^r}{\longrightarrow}QD_r T\zeta.
\edm
The Euler class $e\in\I(M;\zeta )$ is represented by a map $\underline{e}\co M_+\to QT\zeta$.
\begin{cor.}
The following diagram commutes up to homotopy.
\begin{center}
\setlength{\unitlength}{3cm}
\begin{picture}(2.3,1.95)
\put(0.4,0.4){\makebox(0,0){$QD_r T\zeta\times Q(T\zeta\wedge D_{r-1}T\zeta )$}}

 \put(2.6,0.4){\makebox(0,0){$QD_r T\zeta$}}
\put(0.4,1.8){\makebox(0,0){$M_+$}}
\put(2.6,1.8){\makebox(0,0){$N_+$}}
\put(0.6,1.8){\vector(1,0){1.8}}
\put(0.4,1.65){\vector(0,-1){0.4}}
\put(2.6,1.65){\vector(0,-1){0.4}}
\put(1.5,1.83){\makebox(0,0)[b]{$f_+$}}
\put(1.5,0){\makebox(0,0){$QD_r T\zeta\times QD_r T\zeta$}}
\put(0.5,0.3){\vector(1,-1){0.25}}
\put(2.25,0.05){\vector(1,1){0.25}}
\put(0.4,1.1){\makebox(0,0){$QD_r T\zeta\times QT\zeta\times QD_{r-1}T\zeta$}}
\put(2.6,1.1){\makebox(0,0){$QT\zeta$}}
\put(0.4,0.95){\vector(0,-1){0.4}}
\put(2.6,0.95){\vector(0,-1){0.4}}
\put(0.45,0.75){\makebox(0,0)[l]{$1\times\eta$}} \put(2.65,0.75){\makebox(0,0)[l]{$h^r$}}
\put(0.45,1.45){\makebox(0,0)[l]{$\bigl(\underline{\mu_{r+1}(f)}\times\underline{e}\times\underline{\mu_r(f)}\bigr) $}}
\put(2.65,1.45){\makebox(0,0)[l]{$\underline{f}$}}
\end{picture}
\end{center}
\end{cor.}
Here $\eta$ comes from the smash product pairing $QX\times QY\to Q(X\wedge Y)$. The unlabelled map on the left is $1\times Q\xi$, where $\xi \co T\zeta\wedge D_{r-1}T\zeta\to D_r T\zeta$ is a standard map (see the proof of Proposition \ref{properties}). The unlabelled map on the right is the loop product, which corresponds to addition in the bordism group.

\end{document}